\documentclass[12pt]{amsart}
\usepackage{amssymb,amsfonts,amsthm}

\newcommand{\dist}{{\mathrm{dist}}}
\newcommand{\bi}{\bibitem}
\newcommand{\nb}{\newblock}

\newcommand{\la}{\langle\,}
\newcommand{\ra}{\,\rangle}
\newcommand{\zz}{{\mathbb{Z}}}
\newcommand{\be}[1]{\begin{equation}\label{#1}}
\newcommand{\ee}{\end{equation}}
\newcommand{\pp}{{\cal P}}
\newcommand{\rr}{{\cal R}}
\newcommand{\dd}{{\cal D}}

\newcommand{\tr}{\cal T}
\newcommand{\ep}{\epsilon}
\newcommand{\ve}{\varepsilon}
\newcommand{\topp}[1]{\left\lceil{#1}\right\rceil}
\newcommand{\bott}[1]{\left\lfloor{#1}\right\rfloor}
\newcommand{\supp}{\mathop{\mathbf{supp}}}

\newtheorem{thm}{Theorem}[section]
\newtheorem{lm}[thm]{Lemma}

\newtheorem{cy}[thm]{Corollary}

\theoremstyle{definition}
\newtheorem{defn}[thm]{Definition}

\newtheorem{quest}[thm]{Question}

\newtheoremstyle{citing}
  {3pt}
  {3pt}
  {\itshape}
  {}
  {\bfseries}
  {.}
  {.5em}
  {\thmnote{#3}}

\theoremstyle{citing}

\theoremstyle{remark}
\newtheorem{example}{Example}

\let\cal\mathcal

\begin{document}

\title[Metric on diagram groups and uniform embeddings]
{Metrics on diagram groups and uniform embeddings in a Hilbert space}
\author{G.~N.~Arzhantseva}
\address{Section de Math\'ematiques\\
    Universit\'e de Gen\`eve\\ CP 64, 1211 Gen\`eve 4\\
    Switzerland}
\email{Goulnara.Arjantseva@math.unige.ch}
\author{V.~S.~Guba}
\address{Department of mathematics\\
    Vologda State University\\ 6 S. Orlov St.\\ Vologda\\ 160600\\
    Russia}
\email{guba@uni-vologda.ac.ru}
\author{M.~V.~Sapir}
\address{Department of mathematics\\
    Vanderbilt University\\ USA}
\email{m.sapir@vanderbilt.edu}

\thanks{The work of the first two authors has been supported by the Swiss National Science
 Foundation, No.~PP002-68627. They also acknowledge the hospitality of
the Centre de Recerca Matem\`atica (Barcelona) during autumn 2004. The work of the third author has been supported in part by
 the NSF grant DMS 0245600.}


\begin{abstract}
We give first examples of finitely generated groups having an
intermediate, with values in $(0,1)$, Hilbert space compression
(which is a numerical parameter measuring the distortion required to
embed a metric space into  Hilbert space). These groups include certain
diagram groups. In particular, we show that the Hilbert space
compression of Richard Thompson's group $F$ is equal to $1/2$, the
Hilbert space compression of $\zz\wr\zz$ is between $1/2$ and $3/4$,
and the Hilbert space compression of $\zz\wr(\zz\wr\zz)$ is between
0 and $1/2$. In general, we find a relationship between the growth of
$H$ and the Hilbert space compression of $\zz\wr H$.
\end{abstract}

\maketitle

\section{Introduction}

The study of uniform embeddings of metric spaces into Hilbert space was
initiated by Gromov.

\begin{defn}\label{def1}
Let $(\Gamma,d)$ be a metric space. Let $\mathcal{H}$ be a separable
Hilbert space. A map $f\colon\Gamma\to\mathcal{H}$ is said to be a
{\em uniform embedding\/}~{\rm\cite{Gr}} if there exist non-decreasing
functions $\rho_1$, $\rho_2$ from $[0,\infty)$ to itself such that
\begin{itemize}
\item[(1)] $\rho_1(d(x,y))\le \|f(x)-f(y)\|_{\mathcal{H}}\le \rho_2(d(x,y))$
for all $x,y\in\Gamma$;
\item[(2)] $\lim_{r\to+\infty}\rho_i(r)=+\infty$ for $i=1,2$.
\end{itemize}
\end{defn}

Gromov suggested that a finitely generated group $G$ (viewed as a metric space with a
word length metric) uniformly embeddable into Hilbert space should satisfy the Novikov
Conjecture~\cite[page~67]{Gr1}. That was proved in \cite{Yu}, see also~\cite{skandalis,Hig}.


Guoliang Yu introduced a F\o lner-type condition on finitely generated groups $G$,
called property A, which is a weak form of amenability and which guaranties the
existence of a uniform embedding of the metric space into Hilbert space~\cite{Yu}.
That property is interesting in itself because it turned out to be equivalent to
the exactness of $G$, and to the existence of an amenable action of $G$ on a compact
space (see \cite{Tu} for a detailed discussion of property A).

Among ``classical"  groups for which property A has not been proved so far, is the
R.\,Thompson group $F$, which is the group of all piecewise linear orientation
preserving self-homeomorphisms of the unit interval with finitely many dyadic
singularities and all slopes integer powers of $2$.

It is not known whether $F$ is amenable or not, so the question about a weak
amenability is interesting.

It is worth noticing that, by a result of Niblo and Reeves, a group acting properly
and cellularly by isometries on a CAT(0) cubical complex has the Haagerup
property~\cite{Niblo}. In particular, by a result of Farley \cite{Farley}, this holds
for diagram groups over finite semigroup presentations. Moreover, this holds for all
countable diagram groups because they can be embedded into a certain ``universal"
diagram group~\cite{GS3} . It follows immediately from the definition of the Haagerup
property that such a group can be embedded uniformly into Hilbert space. However,
the interaction between the Haagerup property and Guoliang Yu's property A is rather
complicated~\cite[Ch.~1.3]{cherix},~\cite{GuKa,skandalis}. Notice that amenable groups
satisfy both the Haagerup property and property A.

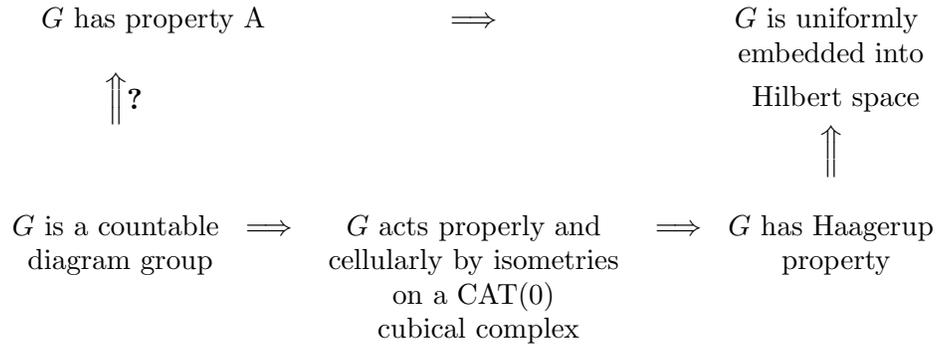
\small
$$
\begin{array}{ccccc}
G\ \hbox{has property A}\hspace{-1cm} & &\Longrightarrow & &G\ \hbox{{\small is
uniformly }}\\
                         &  &      &        & \strut \ \hbox{{\small  embedded
into }}\\
\Big\Uparrow\vcenter{\rlap{{\bf ?}}} &  &      &        &\strut  \ \hbox{{\small
Hilbert space}}\\
& & &    & \Big\Uparrow\\
& & &    & \strut \\
G\ \hbox{is a countable} & \Longrightarrow& G\ \hbox{{\small acts properly and}} &
\Longrightarrow & G \ \hbox{has Haagerup}\\
\ \hbox{diagram group}& &\strut  \ \hbox{cellularly by isometries }& &\strut\
\hbox{property}\\
 & &\strut \ \hbox{on a CAT(0) }  & &\\
 & &\strut \ \hbox{cubical complex}  & &
\end{array}
$$
\normalsize
\bigskip

Guentner and Kaminker \cite{GuKa} introduced a natural quasi-isometry invariant of a
group that shows how close to a quasi-isometry can a uniform embedding of a group
into a Hilbert space be.

\begin{defn}{\rm (cf.\cite[Def.~2.2]{GuKa})}
The Hilbert space compression of a finitely generated discrete group $G$ is the number
$R(G)$, which is the supremum of all $\alpha\ge0$ for which there exists a uniform
embedding of $G$ into a Hilbert space with $\rho_1(n)=Cn^\alpha$ with a constant $C>0$
and linear $\rho_2$ (see Definition \ref{def1}).
\end{defn}

By ~\cite{GuKa}, Hilbert space compression strictly greater than 1/2 implies property A.
Any group that is not uniformly embeddable into a Hilbert space (such groups exist
by \cite{random}) has Hilbert space compression 0. It is proved in \cite{GuKa} that
free groups have Hilbert space compression 1 (although there are no quasi-isometric
embeddings of a free group of rank $>1$ into a Hilbert
space~\cite{bourg,linial}\footnote{We thank D.~Sonkin for bringing these papers to our
attention.}). More generally, by Campbell and Niblo \cite{nibloS}, any discrete group
acting properly, co-compactly on a finite dimensional CAT(0) cubical complex has Hilbert
space compression 1.

Note that till now there were no examples of groups with compression
strictly between 0 and 1.

One of our main results is the following.

\begin{thm}\label{Fhalf}
The Hilbert space compression of R.\,Thompson's group $F$ equals $1/2$.
\end{thm}

In particular, this shows that the  result by Guentner and Kaminker
on groups with compression strictly greater than $1/2$ cannot be
applied to $F$. However, one can extract a stronger fact
from~\cite{GuKa}. Suppose that a finitely generated group $G$ has a
uniform embedding into Hilbert space with linear $\rho_2(n)$ and
$\rho_1\gg\sqrt{n}$ (i.e.
$\lim_{n\to\infty}\rho_1(n)/\sqrt{n}=\infty$). The proof
of~\cite[Theorem~3.2]{GuKa} shows that property A holds for $G$ in
this case. For R.\,Thompson's group $F$ we actually show that
$\rho_1\le C\sqrt{n} \log n$ for some constant $C$.

\begin{quest}
\label{sqd}
Is it true that there exists a uniform embedding of $F$ into Hilbert space
with linear $\rho_2$ and $\rho_1\gg\sqrt{n}$?
\end{quest}

If the answer is positive, then $F$ has property A.
\vspace{1ex}


It is known \cite{GS1,GS5,Farley} that geometric properties of $F$
are better understood when $F$ is considered as a diagram group.
Recall that a diagram group is the fundamental group of the space of
positive paths on a directed 2-complex~\cite{GS3} (an equivalent
definition in terms of semigroup presentations from~\cite{GS1} is
given below). For example, $F$ is the diagram group of the dunce hat
considered as a directed $2$-complex with one vertex, one edge $x$ and
one directed cell $x=x^2$.

Our interest in uniform embeddings for diagram groups was stimulated
in part by a connection with certain metric properties of diagram
groups which have been studied independently before. Recall that
elements of a diagram group $G$ can be represented by diagrams, which are
plane cellular complexes subdivided into a number of regions (cells).
This allows one to introduce a canonical {\em diagram metric\/} $\dist_d$ on
$G$ such that the $\dist_d(g_1,g_2)$ is the number of cells in the diagram
with minimal number of cells representing the element $g_1^{-1}g_2$. This
metric is left invariant. It is proper provided the corresponding directed
2-complex (semigroup presentation) is finite or when the diagram group is
finitely generated. It is known (this is proved in Burillo~\cite{burillo} using
a different terminology) that for the R.\,Thompson group $F$ considered as a
diagram group of the dunce hat, the diagram metric is bi-Lipschitz equivalent
to the word metric. We say that a finitely generated diagram group satisfies
{\em property B\/} if its diagram metric is bi-Lipschitz equivalent to the word
metric.

\begin{thm}\label{dghalf}]
The Hilbert space compression of a finitely generated diagram group with property B
is at least $1/2$.
\end{thm}

\begin{quest} \label{prb} Do all finitely generated diagram groups satisfy property B?
\end{quest}

We think that is a very interesting question.
\vspace{2ex}

In \cite{GS1,GS3}, it is proved that there are (finitely presented) diagram groups
containing all countable diagram groups as subgroups. Such diagram groups are called
{\em universal\/}. One of these diagram groups (in this paper it is denoted by $U$)
corresponds to the semigroup presentation $\la x,a\mid  x^2=x^3,ax=a\ra$.

\begin{thm}\label{BurU}
The universal diagram group $U$ satisfies property B.
\end{thm}

Notice that Theorem \ref{BurU} does not imply (at least directly) a
positive answer to Question \ref{prb}. Indeed, $U$ contains the direct
product $G={\mathbb F}_2\times{\mathbb F}_2$ of two free groups of rank 2
(which is a countable diagram group). By a result of Mikhailova \cite{Mikh},
the group $G$ contains a non-recursively distorted subgroup $H$. Hence the word
metric on a diagram group may not be equivalent to the metric induced by the
embedding of $H$ into $U$. Certainly Theorem \ref{BurU} implies that all
undistorted subgroups of $U$ satisfy property B. Clearly the diagram metric
on a finitely generated diagram group does not exceed a constant times the word
metric. An easy argument shows that for any finitely generated diagram group,
there exists a recursive function $f(n)$ such that the word metric does not
exceed $f(\dist_d)$.

Besides the ~R.\,Thompson group $F$, the restricted wreath product
$\zz\wr\zz$ is another typical representative of the class of diagram groups.
It corresponds to the presentation
$\la a,b,b_1,b_2,c\mid ab=a,bc=c,b=b_1,b_1=b_2, b_2=b\ra$ \cite{GS5}.
It also satisfies property B. Hence Theorem \ref{dghalf} can be applied. The
group $\zz\wr\zz$ is amenable, so it satisfies property A, but the
next theorem shows that the compression of that group is not $1$.

\begin{thm} \label{wr}
The Hilbert space compression of the restricted wreath product
$\zz\wr\zz$ belongs to $[1/2,3/4]$.
\end{thm}

\begin{quest} What is the Hilbert space compression of $\zz\wr\zz$?
\end{quest}

Finally we prove the following general result.

\begin{thm} \label{wr1} Let $H$ be a finitely generated group with
a su\-per-po\-ly\-no\-mial growth function. Then the Hilbert space compression of
$\zz\wr H$ is at most $1/2$.
\end{thm}

This theorem immediately implies

\begin{cy} The Hilbert space compression of the group $\zz\wr(\zz\wr\zz)$
is at most $1/2$.
\end{cy}

Note that the group $\zz\wr(\zz\wr\zz)$ is amenable and so it satisfies property A.
Thus property A does not imply that Hilbert space compression is bigger than $1/2$.
It is not known whether $\zz\wr(\zz\wr\zz)$ is a diagram group (the most probable answer
seems to be negative). Thus Theorem \ref{dghalf} does not apply and we do not know if the
Hilbert space compression of this group is exactly $1/2$.

\begin{quest} What is the Hilbert space compression of the group $\zz\wr(\zz\wr\zz)$?
Are there amenable groups with Hilbert space compression 0 or arbitrary close to 0?
\end{quest}


Finally, we would like to mention that our proof of Theorem~\ref{dghalf} yields actually a
stronger result. Namely, a lower bound on the the {\it equivariant\/} Hilbert space
compression $R_G(G)$ which is the Hilbert space compression defined by restricting
to $G$-equivariant uniform embeddings of $G$ into Hilbert spaces equiped with actions
of $G$ by {\it affine isometries\/}~\cite[Section 5]{GuKa}.

\begin{thm}\label{B-equiv}
The equivariant Hilbert space compression of a finitely generated diagram group
with property B is at least $1/2$.
\end{thm}

The estimate $1/2$ cannot be improved since $R_G(G)>1/2$ implies amenability~\cite{GuKa}
and there are non-amenable diagram groups with property B~\footnote{We thank Y. de Cornulier
for this remark.}. Moreover, one can see from the proof that the diagram metric of a
finitely generated diagram group is a conditionally negative definite function on the group.
Recall that a function $\psi : G\to \mathbb{R}$, satisfying $\psi(g)=\psi(g^{-1})$
is said to be conditionally negative definite if $\sum\psi(g_ig_j^{-1})c_ic_j\le 0$
for all $n\in\mathbb{N}$ and $g_1,\dots ,g_n\in G, c_1, \dots , c_n\in \mathbb{R}$ such that
$\sum c_i=0$. The existence of a continious conditionally negative definite function $\psi : G\to [0,\infty)$ with
$\lim_{g\to\infty}\psi(g)=\infty$ is equivalent to the a-T-menability ~\cite{ch-b-v}, hence
to the Haagerup approximation property~\cite{cherix}. For more applications see~\cite{GuK} and references therein.


\section{Preliminary information about diagram groups}

\subsection{Diagram groups}\label{diagram}

Let us briefly recall the concept of a diagram group and some
terminology from this area. Details can be found in~\cite{GS1}. For
the reasons of the present paper, it is enough to use the definition
in terms of semigroup presentations given in \cite{GS1} rather than
an equivalent definition from \cite{GS3}.

Let $X$ be an alphabet. We denote by $X^*$  the set of all words over $X$
and by $X^+$ the set of all nonempty words. The length of a word $w$
is denoted by $|w|$.

Let $\pp=\la X\mid\rr\ra$ be a semigroup presentation. Here $\rr$ consists of
ordered pairs of nonempty words over $X$. We always assume that if $(u,v)\in\rr$,
then $(v,u)\notin\rr$.

A triple of the form $(p,u=v,q)$ is called an {\em atomic $2$-path\/}, where
$p,q\in X^*$, $(u,v)\in\rr\cup\rr^{-1}$. The following labelled plane graph is
called an {\em atomic diagram\/} (over $\pp$):

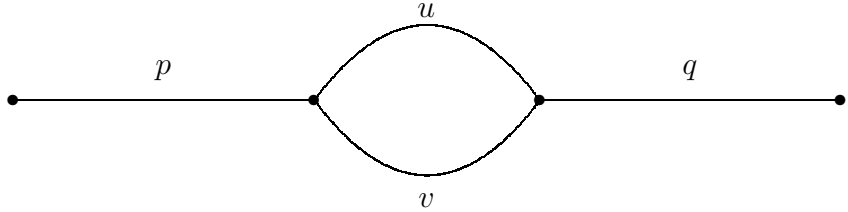
\begin{center}
\unitlength=1.00mm
\special{em:linewidth 0.4pt}
\linethickness{0.4pt}
\begin{picture}(110.66,26.33)
\put(0.00,13.33){\circle*{1.33}}
\put(40.00,13.33){\circle*{1.33}}
\put(70.00,13.33){\circle*{1.33}}
\put(110.00,13.33){\circle*{1.33}}
\bezier{200}(40.00,13.33)(55.00,33.33)(70.00,13.33)
\bezier{200}(40.00,13.33)(55.00,-6.67)(70.00,13.33)
\put(20.00,17.33){\makebox(0,0)[cc]{$p$}}
\put(55.00,25.33){\makebox(0,0)[cc]{$u$}}
\put(90.00,17.33){\makebox(0,0)[cc]{$q$}}
\put(55.00,0.00){\makebox(0,0)[cc]{$v$}}
\put(0.33,13.33){\line(1,0){39.67}}
\put(70.00,13.33){\line(1,0){40.00}}
\end{picture}
\end{center}

\noindent Here each segment labelled by a word $w$ is subdivided into $|w|$ edges
labelled by the letters of $w$. (Segments have a left-to-right orientation.)

Notice that each atomic diagram has the {\em top path\/} (in our case it is
labelled by $puq$) and the {\em bottom path\/} (labelled by $pvq$). This atomic
diagram is defined uniquely by a $2$-path $(p,u=v,q)$ so we will often use
this notation for the atomic diagram itself.

Suppose that we have a sequence $\Delta_1$, \dots, $\Delta_k$ of atomic
diagrams. Let us consider the case when the bottom label of each $\Delta_i$
($1\le i<k$) coincides with the top label of $\Delta_{i+1}$. We can thus
identify the bottom path of $\Delta_i$ with the top path of $\Delta_{i+1}$
($1\le i<k$) to obtain the new labelled plane graph denoted by
$\Delta:=\Delta_1\circ\cdots\circ\Delta_k$ (the {\em concatenation of\/}
the above {\em atomic diagrams\/}). One can naturally define the top path
(which will be the top path of $\Delta_1$) and the bottom path (the bottom
path of $\Delta_k$) of $\Delta$.

For every word $w\in X^+$, we define the {\em trivial diagram\/} denoted by
$\ve(w)$, which is the interval labelled by $w$. Its top and bottom paths
coincide with itself.

\begin{defn}
A {\em diagram $\Delta$ over a semigroup presentation $\pp$\/} is either a trivial
diagram, or a concatenation of atomic diagrams over $\pp$. The diagrams are
considered identical whenever they are isotopic as labelled plane graphs.
\end{defn}

For a diagram $\Delta$, its top and bottom paths will be denoted by
$\topp{\Delta}$ and $\bott{\Delta}$, respectively. They start at the
same vertex denoted by $\iota(\Delta)$ and end at the same vertex
denoted by $\tau(\Delta)$. The diagram itself is situated ``between"
its top and bottom paths. If $w'$ is the label of the top path of
$\Delta$ and $w''$ is the label of its bottom path, then we say that
$\Delta$ is a $(w',w'')$-diagram.

Given two diagrams $\Delta'$ and $\Delta''$, one can define their
{\em concatenation\/} denoted by $\Delta'\circ\Delta''$, provided the
labels of $\bott{\Delta'}$ and $\topp{\Delta''}$ coincide.

For every word $w\in X^+$, the set of all $(w,w)$-diagrams over $\pp$
forms a monoid with $\ve(w)$ as a unit.

For any diagram $\Delta$, by $\Delta^{-1}$ we denote the mirror image of
$\Delta$ with respect to the axis connecting $\iota(\Delta)$ and $\tau(\Delta)$).
Clearly, the top (bottom) path of $\Delta$ becomes the bottom (top) path of
$\Delta^{-1}$.
\vspace{1ex}

Here is an example of two diagrams over $\la x\mid x^3=x^2\ra$:

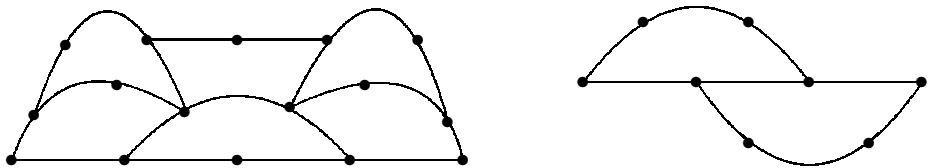
\begin{center}
\unitlength=1.00mm
\special{em:linewidth 0.4pt}
\linethickness{0.4pt}
\begin{picture}(126.33,25.67)
\put(80.33,12.00){\circle*{1.33}}
\put(95.33,12.00){\circle*{1.33}}
\put(110.33,12.00){\circle*{1.33}}
\put(125.33,12.00){\circle*{1.33}}
\bezier{200}(80.33,12.00)(95.33,32.00)(110.33,12.00)
\bezier{212}(95.33,12.00)(110.33,-10.00)(125.33,12.00)
\put(88.33,20.00){\circle*{1.33}}
\put(102.33,20.00){\circle*{1.33}}
\put(102.33,4.00){\circle*{1.33}}
\put(118.33,4.00){\circle*{1.33}}
\put(4.33,1.67){\circle*{1.33}}
\put(19.33,1.67){\circle*{1.33}}
\put(34.33,1.67){\circle*{1.33}}
\put(49.33,1.67){\circle*{1.33}}
\put(64.33,1.67){\circle*{1.33}}
\bezier{180}(19.33,1.67)(34.33,18.67)(49.33,1.67)
\put(27.33,8.00){\circle*{1.33}}
\put(41.33,8.67){\circle*{1.33}}
\bezier{152}(27.33,8.00)(10.33,18.67)(4.33,1.67)
\bezier{152}(41.33,8.67)(60.33,17.67)(64.33,1.67)
\put(7.33,7.67){\circle*{1.33}}
\put(62.33,6.67){\circle*{1.33}}
\put(18.33,11.67){\circle*{1.33}}
\put(51.33,11.67){\circle*{1.33}}
\bezier{228}(7.33,7.67)(16.33,34.67)(27.33,8.67)
\bezier{240}(41.33,8.67)(54.33,35.67)(62.33,6.67)
\put(11.50,17.00){\circle*{1.33}}
\put(22.33,17.67){\circle*{1.33}}
\put(46.33,17.67){\circle*{1.33}}
\put(58.33,17.67){\circle*{1.33}}
\put(34.33,17.67){\circle*{1.33}}
\put(4.33,1.67){\line(1,0){60.33}}
\put(22.33,17.67){\line(1,0){24.00}}
\put(80.33,12.00){\line(1,0){45.00}}
\end{picture}
\end{center}

The label of each edge of these diagrams is $x$.

Notice that each cell in a diagram can be considered as a diagram
itself. If the cell corresponds to a relation $u=v$, we shall denote
this diagram by $(u=v)$.

Let $\Delta$ be a diagram over $\pp$. Suppose that it contains a pair
of cells $\pi'$, $\pi''$ such that the bottom path of $\pi'$ coincides
with the top path of $\pi''$ and the top label of $\pi'$ equals the
bottom label of $\pi''$. In this case we say that the cells $\pi'$
and $\pi''$ form a {\em dipole\/} in $\Delta$. One can define the
operation of removing the dipole. Namely, given a dipole formed by
$\pi'$, $\pi''$, we first remove the common boundary of these cells
and then glue $\topp{\pi'}$ with $\bott{\pi''}$ (the paths we glue
have the same label). We get a new diagram over $\pp$ that has fewer
cells. Proceeding in such a way, we get to a diagram that has no dipoles,
called a {\em reduced\/} diagram.

Kilibarda \cite{KilDiss,Kil} proved that the process of cancelling
dipoles in diagrams is confluent, that is, the result does not
depend on the order in which we remove the dipoles. Two diagrams are
called {\em equivalent} if the corresponding reduced diagrams are
the same.

For every word $w\in X^+$, the set of all reduced $(w,w)$-diagrams
with operation ``concatenation followed by removing all dipoles" is
a group ~\cite[Lemma~5.2]{GS1}, denoted by $\dd(\pp,w)$.

\begin{defn}
The group $\dd(\pp,w)$ is called a {\em diagram group\/} over a semigroup
presentation $\pp$ with {\em base\/} $w$.
\end{defn}

\begin{example}{\rm \cite[Ex.~6.4]{GS1}}
The diagram group $\dd(\pp,x)$ over the semigroup presentation
$\pp=\la x\mid x^2=x\ra$ with base $x$ is R.\,Thompson's group $F$
that has the following group presentation:
$$
\la x_0,x_1,x_2,\ldots\mid x_jx_i=x_ix_{j+1}\ (j>i)\ra.
$$
\end{example}

There is one more important binary operation on diagrams. Given two diagrams
$\Delta_1$, $\Delta_2$ over $\pp$, one can identify $\tau(\Delta_1)$ and
$\iota(\Delta_2)$. This gives a new diagram over $\pp$ denoted by
$\Delta_1+\Delta_2$. This operation is associative but not commutative.

\subsection{Universal diagram groups}\label{universal}

In \cite{GS3}, Guba and Sapir found fi\-ni\-tely presented diagram
groups each of which contains every countable diagram group as a
subgroup. Such diagram groups are called {\em universal\/}. For
instance, the diagram group $U=\dd(\pp,a)$, where $\pp=\la x,a\mid
x^3=x^2,ax=a\ra$ is universal. It has the following Thompson-like
group presentation:
$$
\la x_0,x_1,x_2,\ldots\mid x_jx_i=x_ix_{j+1}\ (j-i>1)\ra.
$$
It is easy to see that $U$ can be generated by $x_0$, $x_1$, $x_2$.
Indeed, $x_n=x_2^{x_0^{n-2}}$ for all $n\ge3$ ($a^b=b^{-1}ab$ by
definition). (In fact this group is finitely presented \cite{GS3}.)

We are going to describe explicitly the procedure of expressing an
element of $U$ as a word in generators $x_0$, $x_1$, $x_2$, \dots\,.
Let $g\in U$ be represented by an $(a,a)$-diagram $\Delta$ over $\pp$.
One can decompose $\Delta$ into a product of atomic diagrams. It is easy
to see that the atomic diagrams that can occur have the form
$(1,a=ax,x^m)^{\ep}$, where $m\ge0$, $\ep=\pm1$ or the form
$(ax^k,x^2=x^3,x^m)^{\ep}$, where $k,m\ge0$, $\ep=\pm1$. In the first case
we assign the identity to the atomic diagram and in the second case we assign
to it the element $x_m^{\ep}$. Multiplying these elements, we get an expression
for the element $g\in U$, see~\cite[Section~6]{GS3}.

Given a diagram $\Xi$ over $\la x\mid x^3=x^2\ra$, we can canonically
assign a diagram $\Delta$ over $\la x,a\mid x^3=x^2,ax=a\ra$ as follows.
First we take the sum $\ve(a)+\Xi$ and then add a number of cells of
the form $a=ax$ on the top to make the top label equal $a$. Then we add
a number of cells of the form $ax=a$ on the bottom in order to make the
bottom label also equal $a$. The result is an $(a,a)$-diagram that
represents an element of the group $U$.

On the above picture, where two diagrams over $\la x\mid x^3=x^2\ra$ are
drawn, the described operation leads to the following elements of $U$,
respectively: $x_3x_1^{-1}x_0^{-1}x_3^{-2}x_1^{-1}$ and $x_1^{-1}x_0$.
Notice that the way to decompose a diagram into the product of atomic
factors is not unique in general. For the first case, we always chose the
rightmost cell that can be included into the next atomic factor. This
procedure will be described later in details.


\section{Main Results}

\begin{defn}
Let $M$ be a set and $p,q\colon M\to[0,\infty)$ be arbitrary functions. Then
$p\preceq q$ if there exists a constant $c>0$ such that $p(m)\le cq(m)$ for every
$m\in M$. These functions are equivalent, $p\sim q$, if $p\preceq q$ and $q\preceq p$.
\end{defn}

Let $G$ be a group generated by a set $\Sigma$. For any $g\in G$, we denote by
$\ell_\Sigma(g)$ the length of the shortest word over $\Sigma^{\pm1}$ that
represents $g$. Notice that if $G$ is finitely generated, then all functions of the
form $\ell_\Sigma$ are equivalent for all finite sets of generators. Therefore, the
function $g\mapsto\ell(g)$ is unique up to equivalence. The function $\ell$ on $G$
defines the {\em word length metric\/} on $G$ as follows: the (word) distance
between $g_1,g_2\in G$ is $\ell(g_1^{-1}g_2)$.

We will also omit the subscript $\Sigma$ on $\ell$ if the generating set of $G$ is clear.

Let $G=\dd(\pp,w)$ be a diagram group over a semigroup presentation $\pp$ with base $w$.
Let $\psi\colon H\hookrightarrow G$ be an embedding of a finitely generated group $H$
into $G$.

\begin{defn}
We say that $\psi\colon H\hookrightarrow G$ is a $B$-{\em embedding\/}, whenever
$\#(h)\sim\ell(h)$ for all $h\in H$. Here $\#(h)$ is the number of cells in the
reduced $(w,w)$-diagram representing $\psi(h)$.
\end{defn}

An obvious argument implies that $\#(g)\preceq\ell(g)$ in all cases. Indeed, let $C$
be the maximum number of cells that represent generators of $H$. Then any element
that has length $n$ in these generators can be represented by a diagram over $\pp$
with at most $Cn$ cells.

However, not every embedding of the above form is a $B$-embedding.

\begin{example}
Let us consider the direct product $G={\mathbb F}_2\times{\mathbb F}_2$ of
two free groups of rank $2$. It is known that this is a diagram group, see,
for example, \cite[Section~8]{GS1}. By a well-known result of Mikhailova \cite{Mikh},
there is a finitely generated subgroup $H$ in $G$ with undecidable membership problem.
This implies that no recursive function $f(n)$ can have the property
$\ell(h)\le f(\#(h))$ for all $h\in H$.
\end{example}

The most important case is when $G=\dd(\pp,w)$ is a finitely generated diagram
group and $H=G$ with the identical embedding $\psi$.

\begin{defn}
A finitely generated diagram group $G=\dd(\pp,w)$ over a semigroup
presentation $\pp$ with base $w$ {\em has property\/} B whenever the
identical embedding $G\hookrightarrow G$ is a $B$-embedding.
\end{defn}

In other words, for all $g\in G$, one has $\#(g)\sim\ell(g)$, where $\#(g)$
is the number of cells in the reduced diagram representing $g$.

Notice that we cannot simply say that $G$ has property B itself.
This concept depends on the presentation $\pp$ and base $w$.

R.\,Thompson's group $F$ has property B as a diagram group over $\la x\mid
x^2=x\ra$ with base $x$. This immediately follows from a result of
Burillo \cite{burillo}.

\begin{proof}[Proof of Theorem \ref{dghalf}] Let $G$ be a finitely generated
group that can be B-embedded into a diagram group. We need to show
that the Hilbert space compression of $G$ is at least $1/2$.

Given a diagram group $\dd(\pp,w)$, let us build the following geometric object
associated with this group.
Let us take all reduced diagrams over $\pp$ that have $w$ as a top label. We
identify all top paths of these diagrams. This gives a 1-path $p$
labelled by $w$.

Suppose that there are two reduced $(w,\cdot)$-diagrams $\Delta_1$,
$\Delta_2$ with decompositions of the form
$\Delta_i\equiv\Delta_i'\circ\Delta_i''$, where $i=1,2$ and
$\Delta_1'$, $\Delta_2'$ are isotopic. Then we identify $\Delta_1'$
with $\Delta_2'$ via this isotopy. We do that for all pairs of
diagrams and all decompositions of them with the above property. The
object we get as a result is a directed 2-complex $\tr$, which turns
out to be a rooted 2-tree in the sense of \cite{GS3}. This directed
$2$-complex can be viewed as a semigroup presentation if we assign
different labels to different edges and consider pairs of words
written on the boundaries of the cells as relations. It is proved in
\cite{GS3} that for any path $q$ in $\tr$ with the same endpoints as
$p$ there exists a unique reduced $(p,q)$-diagram over $\tr$.

Let $\mathbf F$ be the set of all geometric 2-cells of $\tr$. By
${\mathbb R}^{\mathbf F}$ we denote the vector space over ${\mathbb R}$
with ${\mathbf F}$ as a basis. Clearly, this vector space is a
subset in Hilbert space ${\cal H}=\ell_2(\mathbf F)$.

Every element $g$ of the diagram group $\dd(\pp,w)$ can be uniquely
represented by a reduced $(w,w)$-diagram $\Delta$. This diagram can
be naturally embedded into $\tr$. (The top of the diagram is identified
with the path $p$ under this embedding.) Let us assign to $g$ a function
$\phi_g$ from $\mathbf F$ to $\mathbb R$, where $\phi_g(f)=1$ if
$f\in\mathbf F$ is contained in the image of $\Delta$ under the
above embedding and $\phi_g(f)=0$ otherwise.

So we have a mapping $\phi\colon\dd(\pp,w)\to{\cal H}$ from the diagram group
$\dd(\pp,w)$ to ${\cal H}$ defined by the rule $g\mapsto\phi_g$.

As we have already mentioned in the introduction, one can define a
canonical {\em diagram metric\/} on $\dd(\pp,w)$ as follows: given
two elements $g_1,g_2\in\dd(\pp,w)$, one can define the {\em diagram
distance\/}, denoted by $\dist_d(g_1,g_2)$ between these elements as
the number of cells in the reduced diagram over $\pp$ representing
$g_1^{-1}g_2$.

Suppose that the diagram distance between two diagrams $\Delta_1$,
$\Delta_2$ from $\dd(\pp,w)$ equals $n$. Let us consider the images
of the diagrams $\Delta_i$ ($i=1,2$) in $\tr$. They can be
decomposed as $\Delta_i\equiv\Psi\circ\bar\Delta_i$ ($i=1,2$), where
$\Psi$ is the ``greatest common divisor" of $\Delta_1$ and
$\Delta_2$. In this case we see that $\bar\Delta_1$, $\bar\Delta_2$
do not have common cells in $\tr$. The total number of cells in
$\bar\Delta_1$ and $\bar\Delta_2$ equals $n$. Since
$\phi_{g_1}-\phi_{g_2}$, as an element of ${\mathbb R}^{\mathbf F}$
is a vector whose coordinates are $0, 1$ or $-1$, and the number of
non-zero coordinates is $n$, we conclude that the norm of
$\phi_{g_1}-\phi_{g_2}$ is $\sqrt{n}$.

Now given a finitely generated group $G$, which is B-embedded into a
diagram group $\dd(\pp,w)$, we see that the word length metric in
$G$ is equivalent to the diagram metric induced on $G$ as a subset
in $\dd(\pp,w)$. Therefore, for some positive constants $C_1$, $C_2$
one has inequalities
$$C_1\sqrt{\dist(g_1,g_2)}\le\|\phi_{g_1}-\phi_{g_2}\|\le
C_2\sqrt{\dist(g_1,g_2)}\le C_2\dist(g_1,g_2).$$ Hence the Hilbert
space compression of $G$ is at least $\frac12$. 
\end{proof}

 We will use
the following result which is a generalization of the well-known
parallelogram theorem to higher dimensions \cite{DeMa}. Namely, let
$E^n\subset {\mathbb R}^n$ be an $n$-dimensional hypercube. Suppose
that we have a mapping of the set of vertices of $E^n$ into a metric
space $M$. In this case we will say that we have a {\em skew cube\ }
in $M$. For every edge of $E^n$ (there are exactly $2^{n-1}n$ of
them), by an edge of the skew cube we will mean the distance in $M$
between the images of the endpoints of the edge. Similarly, for each
(long) diagonal of $E^n$ (which connects opposite vertices of $E^n$)
we consider the corresponding {\em diagonal\/} of the skew cube.

\begin{lm}{\rm (Skew Cube Inequality \cite{DeMa})}\label{cube}
For every skew cube in a Hilbert space, the sum of squares of its
diagonals does not exceed the sum of squares of its edges.
\end{lm}

\begin{proof}[Proof of Theorem~\ref{Fhalf}] We need to prove that the
Hilbert space compression of R.\,Thompson's group $F$ equals $1/2$.

The fact that the compression is at least $1/2$, follows from
Theorem \ref{dghalf}.

For any $n\ge0$, let us define $2^n$ elements of $F$ that commute
pairwise. All these elements will be reduced $(x,x)$-diagrams over
$\pp=\la x\mid x^2=x\ra$. For $n=0$, let $\Delta$ be the diagram
that corresponds to the generator $x_0$. Namely, if $\pi$ is
the diagram that consists of one cell of the form $x=x^2$, then
$\Delta$ is $\pi\circ(\ve(x)+\pi)\circ(\pi+\ve(x))^{-1}\circ\pi^{-1}$ by
definition. It has 4 cells.

Suppose that $n\ge1$ and we have already constructed diagrams $\Delta_i$
($1\le i\le2^{n-1}$) that commute pairwise. For every $i$ we consider two
$(x^2,x^2)$-diagrams: $\ve(x)+\Delta_i$ and $\Delta_i+\ve(x)$. We get $2^n$
spherical diagrams with base $x^2$ that obviously commute pairwise. It remains
to conjugate them to obtain $2^n$ spherical diagrams with base $x$ having the
same property. Namely, we take $\pi\circ(\ve(x)+\Delta _i)\circ {\pi}^{-1}$ and
$\pi\circ(\Delta_i+\ve(x))\circ\pi ^{-1}$.

Let us denote the elements of $F$ obtained in this way by $g_i$ ($1\le i\le2^n$).
These elements define a $2^n$-dimensional skew cube in $F$. It follows easily from
the construction that each $g_i$ has exactly $2n+4$ cells as a diagram. So the word
length of each $g_i$ is $O(n)$. Now for any $\ep_i=\pm1$ ($1\le i\le2^n$) we consider
the product of the form $g=g_1^{\ep_1}\cdots g_{2^n}^{\ep_{2^n}}$. It is easy to see
from definitions that the diagram that represents $g$ has the form
$\Gamma _n\circ(\Delta^{\ep _1}+\cdots+\Delta^{\ep_{2^n}})\circ\Gamma _n^{-1}$,
where $\Gamma_n$ is defined by induction in the following way: $\Gamma _0=\pi$,
$\Gamma_{k+1}=\pi\circ(\Gamma_k+\Gamma_k)$ ($k\ge0$). In particular,
the number of cells in $\Gamma_n$ equals $2^n-1$ and so $g$ is represented by a
diagram with exactly $2(2^n-1)+4\cdot2^n=3\cdot2^{n+1}-2$ cells. Since $F$ satisfies
property $B$, the word length of an element $g=g_1^{\ep_1}\cdots g_{2^n}^{\ep_{2^n}}$
will be at least $C2^n$, where $C>0$ is a constant that does not depend on the
$\ep_i$'s.

Now consider a uniform embedding of $F$ into a Hilbert space ${\cal
H}$ with linear $\rho_2$. In the image of our skew cube in $F$
formed by $g_1$, \dots, $g_{2^n}$, each edge will be equal to
$O(n)$. The Skew Cube Inequality implies that there exists a
diagonal of the corresponding skew cube in ${\cal H}$ that does not
exceed $O(n)\cdot\sqrt{2^n}=O(n2^{\frac{n}{2}})$. This means that
some points in $F$ that were at distance $d\ge C2^n$ from each other
will be mapped to points in ${\cal H}$ at distance
$O(n2^{\frac{n}{2}})=O(\sqrt{d}\log_2 d)$. Therefore, the
compression cannot exceed $1/2$. \end{proof}

In order to prove Theorem \ref{BurU}, we need the  following simple
construction, called the {\em rightmost decomposition\/} of a
diagram. The idea of such a decomposition applied to the
presentation $\la x\mid x^2=x\ra$ was used in \cite{GS2} to get a
new normal form for the elements of R.\,Thompson's group $F$.

\begin{defn}{\rm (Rightmost decomposition of a diagram)}
Let $\Delta$ be a diagram over a semigroup presentation $\pp$. A
{\em rightmost decomposition\/} of the diagram $\Delta$ into a
product of atomic factors is defined as follows. We proceed by
induction on the number of cells in $\Delta$. If $\Delta$ has no
cells, then the product has no factors. Let $\Delta$ have cells.
Then it has at least one top cell (a cell whose top path is a
subpath in the top path of $\Delta$). Let us consider the rightmost
top cell $\pi$ of $\Delta$. Let $\topp{\Delta}=p\topp{\pi}q$. Then
$\Delta$ is a concatenation of an atomic diagram
$\Pi=\ve(p)+\pi+\ve(q)$ and some diagram $\Delta'$ that has fewer
cells. (Here $\Delta'$ is obtained from $\Delta$ by deleting the
path $\topp{\pi}$.) Taking the product of $\Pi$ and the rightmost
decomposition of $\Delta'$, we obtain the rightmost decomposition of
$\Delta$.
\end{defn}

\begin{proof}[Proof of Theorem~\ref{BurU}] Let $\pp=\la x,a\mid
x^3=x^2,ax=a\ra$. We need to prove that the universal diagram group
$U=\dd(\pp,a)$ has property B.

Let $\Delta$ be an $(a,a)$-diagram over $\pp$ with $N$ cells. It
suffices to prove that the element $g\in U$ represented by $\Delta$
has length at most $KN$ in generators $x_0$, $x_1$, $x_2$, where
$K>0$ is a constant independent on $g$.

Let $e_0$ be the top edge of $\Delta$. Suppose that there is an
$(a,ax)$-cell whose top edge coincides with $e_0$. Then we denote by
$e_1$ the edge labelled by $a$ on the bottom path of this cell. If
$e_1$ is the top path of an $(a,ax)$-cell, then $e_2$ denotes the
edge on the bottom path of this cell labelled by $a$. Proceeding in
such a way, we finally obtain a sequence of edges $e_0$, \dots,
$e_k$ ($k\ge0$).

Analogously, changing in the previous paragraph top by bottom, we define the sequence
of edges $f_0$, \dots, $f_m$  ($m\ge0$), where $f_0$ is the bottom edge of $\Delta$.
A very easy geometric observation is that $e_k$ must coincide with $f_m$. It is also
easy to see that $\Delta$ is a concatenation of the form
$\Delta=\Delta_1\circ(\ve(a)+\Delta')\circ\Delta_2$, where $\Delta_1$ consists of
$(a,ax)$-cells only, $\Delta_2$ consists of $(ax,a)$-cells only and $\Delta'$ is an
$(x^k,x^m)$-diagram over $\pp'=\la x\mid x^3=x^2\ra$.

Let us apply the rightmost decomposition to $\ve(a)+\Delta'$. Each factor is an
atomic diagram of the form $(ax^s,x^2=x^3,x^t)^\ep$, where $s,t\ge0$, $\ep=\pm1$.
According to the description of $U$ given in subsection~\ref{universal}, this atomic
diagram corresponds to $x_t^\ep$. (Notice also that no generators correspond to atomic
diagram with $(a,ax)$-cells and/or their inverses.)

Therefore, the rightmost decomposition of $\Delta'$ allows us to decompose $g$ as a
product (in $U$) of the form
\be{epss}
g=x_{j_1}^{\ep_1}\cdots x_{j_r}^{\ep_r}\hbox{ with }\ep_i=\pm 1,
\ee
where $r=N-k-m$ is the number of cells in $\Delta'$. We need to establish some easy
properties of the subscripts and the superscripts in (\ref{epss}).

Suppose that (\ref{epss}) contains a subword of the form $x_i^{\ep}x_j$. Then we claim
that $i\le j+1$. Indeed, otherwise the cell that corresponds to $x_j$ is located to
the right of the cell corresponding to $x_i^{\ep}$. (One can also see that $i>j+1$
would imply that $x_i^{\ep}x_j$ equals $x_jx_{i+1}^{\ep}$ in $U$.)

Now suppose that (\ref{epss}) contains a subword of the form
$x_i^{\ep}x_j^{-1}$. In this case we claim that $i\le j+2$. (Now
$i>j+2$ would also imply that the cell that corresponds to $x_j^{-1}$ is located to
the right of the one corresponding to $x_i^{\ep}$. The element $x_i^{\ep}x_j^{-1}$
would be also equal to $x_j^{-1}x_{i-1}^{\ep}$ in $U$.)

For every $j\ge0$ we replace each letter $x_j$ in (\ref{epss}) by
the word $u_j^{-1}v_ju_j$, where $u_j=1$, $v_j=x_j$ if $j=0,1$ and
$u_j=x_0^{j-2}$, $v_j=x_2$ if $j\ge2$. We get the following word
$W$ in generators $x_0$, $x_1$, $x_2$:
$$
W=u_{j_1}^{-1}v_{j_1}^{\ep_1}u_{j_1}u_{j_2}^{-1}v_{j_2}^{\ep_2}\cdots
v_{j_r}^{\ep_r}u_{j_r}.
$$

Since $x_j=u_j^{-1}v_ju_j$ in $U$, the word $W$ represent the same
element $g\in U$. Now let $W_i$ denote the freely irreducible form
of the word $u_{j_i}u_{j_{i+1}}^{-1}$ for all $1\le i<r$. It is
also convenient to define the words $W_0=u_{j_1}^{-1}$ and
$W_r=u_{j_r}$. We are going to estimate the length of the word
$$
\overline{W}=W_0v_{j_1}^{\ep_1}W_1\cdots v_{j_r}^{\ep_r}W_r.
$$

We know that $|W_0|=j_1\le k-2$ (the first cell chosen in
$\Delta'$ is at a distance $j_1$ from the rightmost point of
$\Delta$ and the length of its top path is at least $2$).
Analogously, $|W_r|=j_r\le m-2$. Also $|v_{j_1}|+\cdots+|v_{j_r}|=r$. So
we estimate some part of the length of $\overline{W}$ as follows:
\be{wv}
|W_0|+|W_r|+|v_{j_1}|+\cdots+|v_{j_r}|\le k+m+r-4<N.
\ee

It remains to estimate the sum $|W_1|+\cdots+|W_{r-1}|$. It is
easy to see from the definitions that $|W_i|\le|j_{i+1}-j_i|$ for
all $1\le i<r$. Let $I$ be the set of all $1\le i<r$ such that
$j_i>j_{i+1}$. Then $|W_1|+\cdots+|W_{r-1}|\le\sum_{1\le i<r}|j_{i+1}-j_i|=S_1+S_2$,
where $S_1=\sum_{i\in I}(j_i-j_{i+1})$ and $S_2=\sum_{i\notin I}(j_{i+1}-j_i)$.

For every $i\in I$, we know that $j_i-j_{i+1}\le2$. Therefore,
$S_1\le2|I|\le2(r-1)<2r$. On the other hand,
$$
S_2-S_1=\sum_{1\le i<r}(j_{i+1}-j_i)=j_r-j_1\le j_r\le m-2<k+m=N-r.
$$
This gives $S_2=S_1+(S_2-S_1)<N+r$ and so $S_1+S_2<N+3r\le4N$.
Using (\ref{wv}), we finally have
$$
|\overline{W}|=(|W_0|+|W_r|+|v_{j_1}|+\cdots+|v_{j_r}|)+(|W_1|+\cdots+|W_{r-1}|)<5N.
$$
So our statement is true for $K=5$.

Notice that the constant here is not optimal. \end{proof}

Let us recall the definition of the (restricted) wreath product of groups. Let $G$, $H$
be two groups. If $\phi$ is a function from $H$ to $G$, then the {\em support\/} of $\phi$
is the set $\supp(\phi)=\{\,b\in H\mid \phi(b)\ne1\,\}$. The restricted wreath product
$G\wr H$ of the above groups will consist of all pairs of the form $(b,\phi)$, where $b\in H$
and $\phi$ is a function from $H$ to $G$ with finite support. The group operation is defined
in the standard way, that is, $(b_1,\phi_1)(b_2,\phi_2)=(b_1b_2,\phi_1^{b_2}\phi_2)$, where
$\phi^\beta(b)=\phi(b\beta)$ by definition ($\beta\in H)$. For each of the groups $G$, $H$ we
have a canonical embedding into $G\wr H$. Namely, $a\mapsto(1,\phi_a)$, where $\phi_a(1)=a$,
$\phi_a(b)=1$ whenever $b\ne1$. Also $b\mapsto(b,\mathbf{1})$, where $\mathbf{1}$ is the
function from $H$ to $G$ that takes all elements to 1.

We will deal with restricted wreath products of finitely generated groups. If we fix some
generating sets for $G$ and for $H$, then the union of these sets (via canonical embeddings)
will generate $G\wr H$. Hence we have three length functions: on $G$, on $H$, and on $G\wr H$.
For simplicity, we will always freely speak about lengths in each of these three groups
provided the generating sets for $G$ and for $H$ are clear. The same will concern Cayley
graphs of these groups. The same symbol $\ell$ will be used to denote the length function
for each of these groups. This will not lead to a confusion.

Working with lengths, we need the following fact proved in \cite{Parry} (see Theorem 1.2).

\begin{lm}
\label{Parr}
Let $G$, $H$ be finitely generated groups with fixed generating sets. Then the length
of an element $(b,\phi)\in G\wr H$ is equal to the sum $\sum\ell(\phi(\beta))+M$, where
$\ell(a)$ is the length of $a\in G$ and $M$ is the length of the shortest path in the
Cayley graph of $H$ starting at the identity, ending at $b$ and visiting each element
of $\supp(\phi)$ at least once. The sum is taken over all $\beta\in\supp(\phi)$.
\end{lm}

Let $G=\dd(\pp,w)$ be a diagram group. We say that the base $w$ is {\em protected\/}
if $w$ is a letter such that $w=uwv$ modulo $\pp$ implies that both words $u$, $v$ are
empty (for any $u$, $v$). For any diagram group $G=\dd(\pp,w)$, one can construct a new
semigroup presentation $\pp'$ adding 3 new letters $a$, $b$, $c$ to the alphabet and one
new defining relation $a=bwc$. It is not hard to see that $G$ will be isomorphic to the
diagram group $\dd(\pp',a)$, where the base $a$ is now protected.

It is shown in \cite[Theorem 11]{GS5} that, given a diagram group $G=\dd(\pp,y)$ with
protected base $y$, one can construct a new semigroup presentation $\bar\pp$ adding two
new letters $x$, $z$ to the alphabet and two new defining relations $x=xy$, $yz=z$
(see also Example 10 before the theorem we refer to). In this situation, the diagram
group $\dd(\bar\pp,xz)$ will be isomorphic to the restricted wreath product $G\wr\zz$.
(The concept of a protected base we gave slightly differs from the original one although
the proof goes without any essential changes.) Now we are going to prove that if $G$ has
property B (as a diagram group over $\pp$) then $G\wr\zz$ also has property B (as a diagram
group over $\bar\pp$).

\begin{lm}
\label{GwrZB}
In the above notation, if $G=\dd(\pp,y)$ has property B, then $G\wr\zz=\dd(\bar\pp,xz)$
also has property B.
\end{lm}

\proof
It suffices to prove that any diagram over $\bar\pp$ with $N$ cells has length $O(N)$ in
the group $G\wr\zz$. Any reduced $(xz,xz)$-diagram over $\bar\pp$ can be obtained as
follows. Let us take finitely many $(y,y)$-diagrams over $\pp$ and let $\Xi$ be their
sum. Now let us add a few cells of the form $x=xy$ and a few cells of the form $z=yz$ to
the diagram $\ve(x)+\Xi+\ve(z)$ on the top in order to make the top label equal $xz$. Then
we add a few cells of the form $xy=x$ and $yz=z$ on the bottom to make the bottom label also equal $xz$.
The result is an $(xz,xz)$-diagram $\Delta$. Let us denote by
$o$ the vertex on the bottom path of $\Delta$ that splits this path into the product of two
edges. Clearly, the vertex $o$ is contained in $\Xi$ and splits it into a sum. Let
$\Xi=\Xi_k+\cdots+\Xi_0+\Xi_{-1}+\cdots+\Xi_{-m}$, where each
$\Xi_j$ ($-m\le j\le k$) is a $(y,y)$-diagram and the vertex $o$ subdivides $\Xi$
into the sum of the two summands. All the $\Xi_j$ that form the first summand satisfy
$j\ge0$, the ones that form the second summand satisfy $j<0$. (Notice that each of the
summands can be empty.)

Let $g_j\in G$ be the element of $\dd(\pp,y)$ represented by $\Xi_j$ ($-m\le j\le k$).
We also define $g_j=1$ for any integer $j$ not in $[-m,k]$. Then, according to
\cite[Example 10]{GS5}, the diagram $\Delta$ represents the element $h$ of  the group
$G\wr\la t\ra$ of the following form: $(t^d,\phi)$, where $d$, $\phi$ are defined as
follows. The exponent $d$ on $t$  is the difference between the number of $(z,yz)$-cells
and the number of $(yz,z)$-cells in $\Delta$. The function $\phi$ from $\zz=\la t\ra$ is
defined by the rule $\phi(t^j)=g_j$ for any integer $j$.

It remains to estimate the length of $h$. We will use Lemma \ref{Parr}. We see that
the length of $h$ equals
\be{sumss}
\sum\limits_{j=-m}^k\ell(g_j)+M,
\ee
where $M$ is the length of the shortest path in the Cayley graph of $\zz=\la t\ra$ that
starts at the identity, ends at $t^d$ and visits all those vertices $t^j$ for which
$g_j\ne1$. It is clear that $M\le2(k+m+1)+|d|$. Obviously, $|d|$ does not exceed the
number of cells in $\Delta$. It is also clear that $k+m+1$ is the number of $(x,xy)$-cells
plus the number of $(z,yz)$-cells in $\Delta$. This implies $M=O(N)$.

Since $G=\dd(\pp,y)$ has property B, there exists a constant $C>0$ such that the
length of each $g\in G$ does not exceed $C\cdot \#(g)$, where $\#(g)$ denotes the
number of cells in the reduced diagram over $\pp$ representing $g$. Therefore, the
first summand in~(\ref{sumss}) does not exceed
$C(\#(g_{-m})+\cdots+\#(g_0)+\cdots+\#(g_k))\le CN$. Thus we finally have that the
length of $h$ in $G\wr\zz$ is $O(N)$.
\endproof

\begin{cy}
\label{zwrzb}
Let $\pp=\la a,b,b_1,b_2,c\mid ab=a,bc=c,b=b_1,b_1=b_2,b_2=b\ra$. Then
the restricted wreath product $\zz\wr\zz$ has property B as a diagram group
$\dd(\pp,ac)$.
\end{cy}

Indeed, $\zz$ is the diagram group over $\la b,b_1,b_2\mid b=b_1,b_1=b_2,b_2=b\ra$
with protected base $b$.
\vspace{2ex}

Let $H$ be a finitely generated group with fixed generating set. By its
{\em growth function\/} $\gamma(n)$ (with respect to this set of generators)
we mean the number of elements in the ball of radius $n$ around the identity
in the Cayley graph of $H$.

Theorems \ref{wr} and \ref{wr1} will follow from the next result.

\begin{thm}\label{wr2}
Let $H$ be a finitely generated group with fixed generating set. Suppose
that its growth function satisfies the condition $\gamma(n)\succeq n^k$ for
some $k>0$. Then the Hilbert space compression $\alpha$ of the group $\zz\wr H$
satisfies the following inequality:
$$
\alpha\le\frac{1+k/2}{1+k}.
$$
\end{thm}

\proof Let $x$ be a generator of $\zz$ in the wreath product $\zz\wr H$. Every
element in $\zz\wr H$ is a pair $(b,\phi)$ where $b\in H$, $\phi$ is a function
$g\mapsto x^{m_g}$ from $H$ to $\zz$ with finite support.

Let $B_n$ be the ball of radius $n$ around $1$ in the Cayley graph of $H$.
Consider the set
$$
X_n=\{\,w_b=b^{-1}x^{2n+1-|b|}b\mid b\in B_n\,\},
$$
where $\ell(b)$ is the length of $b$. Obviously, $w_b$ has the form $(1,\phi_b)$ in
$\zz\wr H$, where the support of $\phi_b$ is $\{\,b\,\}$ and $\phi_b(b)=x^{2n+1-\ell(b)}$.
Notice that $\ell(x^r)=|r|$ in $\la x\ra$ for all $r\in\zz$. Hence by Lemma~\ref{Parr},
the length of $w_b$ in $\zz\wr H$ equals $2n+1-\ell(b)+M$, where $M=2\ell(b)$.
Therefore, $\ell(w_b)=2n+1+\ell(g)$ is always between $2n+1$ and $3n+1$. Clearly, all
elements of $X_n$ pairwise commute. The number of them is exactly $\gamma(n)$.

It is easy to see that for every choice of $\epsilon_b\in\{\,1,-1\,\},b\in B_n$, the
length of the element $w=\prod_{b\in B_n}w_b^{\epsilon_b}$ is greater than $n\gamma(n)$.
Indeed, $w$ has the form $(1,\phi)$ in $\zz\wr H$, where $\supp(\phi)=B_n$ and
$\phi(b)=x^{2n+1-\ell(b)}$ for all $b\in B_n$. Hence the length of $w$ is at least
$$
\sum_{b\in B_n}\ell(\phi(b))=\sum_{b\in B_n}(2n+1-\ell(b))>n\gamma(n)
$$
by Lemma~\ref{Parr}.

Consider the skew cube spanned by $X_n$ in $\zz\wr H$ and its image
in ${\cal H}$ under the uniform embedding $f$. The sides of the skew
cube in $\cal H$ do not exceed $Cn$ for some constant $n$ since our
group is finitely generated. So the sum of squares of sides does not
exceed $C^22^{\gamma(n)-1}\gamma(n)n^2$. On the other hand, the number of
diagonals is $2^{\gamma(n)-1}$. Hence by the Skew Cube Inequality,
there exists a diagonal of the skew cube in $\cal H$ which does not
exceed $Cn\sqrt{\gamma(n)}$. Thus there exist two points in the
Cayley graph of $G$ at distance $d\ge n\gamma(n)$ whose images under
$f$ are at distance at most $Cn\sqrt{\gamma(n)}$. If $\rho_1(d)=Kd^\alpha$
for some constant $K>0$ in the definition of the uniform embedding, then
$Cn\sqrt{\gamma(n)}\ge K(n\gamma(n))^\alpha$. Hence
$$
n^{1-\alpha}\succeq(\gamma(n))^{\alpha-1/2}\succeq n^{k(\alpha-1/2)}.
$$
This implies $1-\alpha\ge k(\alpha-1/2)$ and so $\alpha\le(1+k/2)/(1+k)$.
\endproof

\begin{thm}
\label{wrprod} The Hilbert space compression of the restricted wreath product
$\zz\wr\zz$ belongs to $[1/2;3/4]$.
\end{thm}

\proof The lower bound follows from Theorem \ref{dghalf} and the fact that
$\zz\wr\zz$ satisfies property B as a diagram group $\dd(\pp,ac)$ over
$\pp=\la a,b,b_1,b_2,c\mid ab=a,bc=c,b=b_1,b_1=b_2,b_2=b\ra$
(Corollary~\ref{zwrzb}).

The upper bound follows from Theorem \ref{wr2} with $G=\zz$. Indeed, the growth
function of $\zz$ is linear, that is, $\gamma(n)\succeq n^k$ for $k=1$. Hence
$\alpha\le(1+k/2)/(1+k)=3/4$.
\endproof

Recall that a finitely generated group has a {\em su\-per-po\-ly\-no\-mial\/}
growth whenever its growth function $\gamma(n)$ exceeds $n^k$ for all $k\ge1$.

\begin{thm}
Let $H$ be a finitely generated group with su\-per-po\-ly\-no\-mial growth. Then the
Hilbert  space compression of $\zz\wr H$ is at most $1/2$.
\end{thm}

\proof Indeed, $\gamma(n)\succeq n^k$ implies $\alpha\le(1+k/2)/(1+k)$. The right-hand
side of this inequality approaches $1/2$ as $k\to\infty$.
\endproof

\begin{cy} The Hilbert space compression of the group $\zz\wr(\zz\wr\zz)$
is at most $1/2$.
\end{cy}

\end{document}